\documentclass[12pt,leqno]{amsart}

\usepackage{enumitem}
\usepackage{amsmath,amssymb}
\usepackage{amsthm}   
\usepackage[pdftex,
            colorlinks = true,
            linkcolor = red,
            urlcolor  = blue,
            citecolor = blue,
            anchorcolor = blue]{hyperref}
						
\numberwithin{equation}{section}

\newtheorem{theorem}{Theorem}[section]

{Proposition}[section]
\newtheorem{corollary}
{Corollary}[section]
\newtheorem{remark}{Remark}[section]

\newfont{\bb}{msbm10 at 12pt}

\def\BS{\mathbf{S}}

\def\<{\langle}     
\def\>{\rangle}     
\def\wti{\widetilde}
\def\ov{\overline}

\def\sing{{\rm sing}}

\newcommand{\bal}{\begin{align}}      \newcommand{\eal}{\end{align}}
\newcommand{\ba}{\begin{array}}      \newcommand{\ea}{\end{array}}
\newcommand{\bc}{\begin{center}}     \newcommand{\ec}{\end{center}}
\newcommand{\be}{\begin{enumerate}}  \newcommand{\ee}{\end{enumerate}}
\newcommand{\beq}{\begin{eqnarray}}  \newcommand{\eeq}{\end{eqnarray}}
\newcommand{\beQ}{\begin{eqnarray*}} \newcommand{\eeQ}{\end{eqnarray*}}
\newcommand{\bi}{\begin{itemize}}    \newcommand{\ei}{\end{itemize}}
\newcommand{\bt}{\begin{tabular}}    \newcommand{\et}{\end{tabular}}
\newcommand{\bdm}{\begin{displaymath}} \newcommand{\edm}{\end{displaymath}}

\newcommand{\D}{D\!\!\!\!/\,}
\newcommand{\nbs}{\nabla\!\!\!\!/\,}

\newcommand{\RSB}{\mathbf{S}\!\!\!\!/\,}

\begin{document}

\title[First boundary Dirac eigenvalue and boundary capacity potential]{First boundary Dirac eigenvalue and boundary capacity potential}       

\author{Simon Raulot}
\address[Simon Raulot]{Laboratoire de Math\'ematiques R. Salem
UMR $6085$ CNRS-Universit\'e de Rouen
Avenue de l'Universit\'e, BP.$12$
Technop\^ole du Madrillet
$76801$ Saint-\'Etienne-du-Rouvray, France.}
\email{simon.raulot@univ-rouen.fr}

\begin{abstract}
We derive new lower bounds for the first eigenvalue of the Dirac operator of an oriented hypersurface $\Sigma$ bounding a noncompact domain in a spin asymptotically flat manifold $(M^n,g)$ with nonnegative scalar curvature. These bounds involve the boundary capacity potential and, in some cases, the capacity of $\Sigma$ in $(M^n,g)$ yielding several new geometric inequalities. The proof of our main result relies on an estimate for the first eigenvalue of the Dirac operator of boundaries of compact Riemannian spin manifolds endowed with a singular metric which may have independent interest. 
\end{abstract}

\keywords{}

\subjclass{53C20, 53C27, 83C40, 83C57}

\thanks{}

\date{\today}   

\maketitle 
\pagenumbering{arabic}


\section{Introduction}


The Positive Mass Theorem (PMT) is a famous result in mathematical general relativity which states, in its Riemannian version, that the ADM mass of an $n$-dimensional, $n\geq 3$, complete, asymptotically flat manifold $(M^n,g)$ with nonnegative scalar curvature is nonnegative and it is zero if, and only if, $(M^n,g)$ is isometric to the Euclidean space $(\mathbb{R}^n,\delta)$ (see Appendix \ref{SectionCapacity} for the definitions). It has been first proved by Schoen and Yau \cite{SchoenYau1,SchoenYau2} using the minimal surface technique when $3\leq n\leq 7$ and then by Witten \cite{Witten1} for spin manifolds using the Dirac operator (see also \cite{ParkerTaubes,Bartnik1}). 

Using spinor methods, Herzlich \cite{Herzlich1,Herzlich2} noticed that if $M$ has a compact inner boundary $\Sigma$ with induced metric $\gamma$, some control of the first eigenvalue $\lambda_1(\D_\gamma)$ of $(\Sigma^{n-1},\gamma)$ ensures nonnegativity of the mass. More precisely, he proved the following boundary version of the PMT:
\begin{theorem}$($\cite{Herzlich1,Herzlich2}$)$\label{HerzlichPMT-n}
Let $(M^n,g)$ be an $n$-dimensional, complete, spin asymptotically flat manifold with nonnegative scalar curvature and with a compact inner boundary $\Sigma$ satisfying 
\begin{eqnarray}\label{DiracMeanCur-n}
\lambda_1(\D_\gamma)\geq\frac{1}{2}\max_\Sigma H_g.
\end{eqnarray}
Then the ADM mass $m_{ADM}(M,g)$ of $(M^n,g)$ is nonnegative and if it is zero, $(M^n,g)$ is flat, the mean curvature $H_g$ is constant and (\ref{DiracMeanCur-n}) is an equality. 
\end{theorem}
This result highlights in particular a tight relation between the first eigenvalue $\lambda_1(\D_\gamma)$, the mean curvature $H_g$ of $\Sigma$ and the sign of the ADM mass of $(M^n,g)$. Here the mean curvature for hypersurfaces in asymptotically flat manifolds is computed with respect to the unit normal pointing to infinity.  This version of the PMT was a first step to get a Penrose-like inequality. Recall that the Riemannian Penrose Inequality is a conjecture which 
relates the ADM mass of a complete, asymptotically flat manifold $(M^n,g)$ with nonnegative scalar curvature and the area of its outer minimizing minimal boundary. This inequality is known to be true in dimension three by works of Huisken and Ilmanen \cite{HuiskenIlmanen} for connected boundaries and by Bray \cite{Bray} for the general case (see \cite{BrayLee} for the case $4\leq n \leq 7$). An important tool in Bray's approach is provided by a mass-capacity inequality, recently generalized by Miao and Hirsch \cite{MiaoHirsch}, which states that:
\begin{theorem}$($\cite{Bray,MiaoHirsch}$)$\label{MassCapacity}
Let $(M^n,g)$ be an $n$-dimensional, complete, asymptotically flat manifold with nonnegative scalar curvature and compact inner boundary $\Sigma$. If the boundary capacity potential $\phi:M\rightarrow\mathbb{R}$ satisfies 
\begin{eqnarray}\label{CondMH}
-2\frac{\alpha-1}{\alpha-2}\frac{\partial\phi}{\partial\nu}\geq \frac{n-2}{n-1}H_g
\end{eqnarray}
on $\Sigma$ for some $\alpha\in]0,2[$, then
\begin{eqnarray}\label{MCI}
m_{ADM}(M,g)\geq\alpha\,\mathcal{C}_g(\Sigma,M).
\end{eqnarray}
Moreover, equality holds in (\ref{MCI}) if, and only if, $(M^n,g)$ is isometric to the exterior of a rotationally symmetric sphere in the Riemannian Schwarzschild manifold of corresponding mass. 
\end{theorem} 
The boundary capacity potential of $\Sigma$ in $(M^n,g)$ is the function $\phi:M\rightarrow\mathbb{R}$ which satisfies
\begin{equation}\label{Capacity}
\left\lbrace
\begin{array}{ll}
\Delta_g\phi = 0 & \text{ on } M\\
\phi=1 & \text{ at } \Sigma\\
\phi\rightarrow 0 & \text{ as }x\rightarrow\infty
\end{array}
\right.
\end{equation}
and the boundary capacity of $\Sigma$ in $(M^n,g)$ is 
\begin{eqnarray}\label{DefCap}
\mathcal{C}_g(\Sigma,M):=\frac{1}{(n-2)\omega_{n-1}}\int_M|\nabla^g\phi|^2d\mu_g=-\frac{1}{(n-2)\omega_{n-1}}\int_\Sigma\frac{\partial\phi}{\partial\nu}d\mu_\gamma
\end{eqnarray}
where $\nu$ is the unit normal pointing inside $M$. Here $\Delta_g$ and $\nabla^g$ denote respectively the Laplace operator and the gradient of $(M^n,g)$. A natural question, when comparing the conditions (\ref{DiracMeanCur-n}) and (\ref{CondMH}), is whether there exists a relation between the first boundary Dirac eigenvalue, the mean curvature and the boundary capacity potential in this context. We provide here an affirmative answer regarding this question by showing the following general estimate: 
\begin{theorem}\label{DiracEigenvalueBound}
Let $(M^n,g)$ be an $n$-dimensional, spin, complete, asymptotically flat manifold with nonnegative scalar curvature and compact inner boundary $\Sigma=\coprod_{j=1}^N\Sigma_j$. If 
\begin{eqnarray}\label{PosMC}
-2\frac{n-1}{n-2}\frac{\partial\phi}{\partial\nu}>H_g
\end{eqnarray}
holds on $\Sigma$, then
\begin{eqnarray}\label{DiracCapacity}
\lambda_1(\D_\gamma^{j})\geq\frac{1}{2}\min_{\Sigma_{j}}\Big(-2\frac{n-1}{n-2}\frac{\partial\phi}{\partial\nu}-H_g\Big)
\end{eqnarray}
for all $j=1,\cdots,N$. Moreover, equality holds if, and only if, the Riemannian spin manifold $(M^n,\overline{g})$ with $\overline{g}=\phi^{\frac{4}{n-2}} g$ carries a parallel spinor and $\Sigma$ is connected with positive constant mean curvature. 
\end{theorem}
Here $\Sigma_j$ denotes a connected component of $\Sigma$ for $j=1,\cdots,N$ and $\lambda_1(\D^j_\gamma)$ its corresponding first Dirac eigenvalue. The proof of Theorem \ref{DiracEigenvalueBound} relies on the validity of an eigenvalue estimate by Hijazi, Montiel and Zhang \cite{HijaziMontielZhang1} in the context of singular metrics as defined by Mantoulidis, Miao and Tam in \cite{MantoulidisMiaoTam} and which states that if $\Sigma$ is the mean convex boundary of a {\it compact} Riemannian spin manifold $(\Omega^n,\wti{g})$ with nonnegative scalar curvature then
\begin{eqnarray}\label{HMZ}
\lambda_1(\D_\gamma^j)\geq\frac{1}{2}\min_{\Sigma_j} H_{\wti{g}}.
\end{eqnarray}
The proof of this inequality is given in Section \ref{ExtrinsicLowerBoundSingular}. This result is of independent interest and should have several others applications. 

The inequality (\ref{DiracCapacity}) is sharp since, as computed in Appendix \ref{RiemSchwar}, the exterior region $\mathbb{M}_m^n(r_0)$ of a rotationally symmetric sphere in the Riemannian Schwarzschild manifold $(\mathbb{M}^n_m,g_m)$ satisfies the equality case. In fact, in the $3$-dimensional case and when we restrict the topology of the boundary, one can completely characterize this equality case. 
\begin{corollary}\label{EqualityCaseThree}
Let $(M^3,g)$ be a $3$-dimensional, complete, asymptotically flat manifold with nonnegative scalar curvature and whose inner boundary is a union of $2$-spheres. Then equality holds in (\ref{DiracCapacity}) if, and only if, $(M^3,g)$ is isometric to $(\mathbb{M}_m^n(r_0),g_m)$ for some $r_0$. 
\end{corollary}

Another immediate consequence of Theorem \ref{DiracEigenvalueBound} is a nontrivial estimate for minimal boundaries (or even more generally for boundaries with nonpositive mean curvature), namely:
\begin{corollary}\label{MinimalEstimate}
Let $(M^n,g)$ be an $n$-dimensional, spin, complete, asymptotically flat manifold with nonnegative scalar curvature and minimal inner boundary $\Sigma=\coprod_{j=1}^N\Sigma_j$. Then
\begin{eqnarray}\label{DiracCapacityMinimal}
\lambda_1(\D_\gamma^j)\geq\frac{n-1}{n-2}\min_{\Sigma_j}\big(-\frac{\partial\phi}{\partial\nu}\big)
\end{eqnarray}
for all $j=1,\cdots,N$. Moreover, equality holds if, and only if, the Riemannian spin manifold $(M^n,\overline{g})$ with $\overline{g}=\phi^{\frac{4}{n-2}} g$ carries a parallel spinor and $\Sigma$ is connected with positive constant mean curvature. 
\end{corollary}
A more difficult question is to find a relation between the first eigenvalue of the Dirac operator on $(\Sigma^{n-1},\gamma)$ and the boundary capacity of $\Sigma$ in $(M^n,g)$. When the boundary $\Sigma$ is connected and the normal derivative of the boundary capacity potential is constant along $\Sigma$, Theorem \ref{DiracEigenvalueBound} provides an answer.
\begin{corollary}\label{DiracCapacityConstante}
Let $(M^n,g)$ be an $n$-dimensional, spin, complete, asymptotically flat manifold with nonnegative scalar curvature and connected, compact inner boundary $\Sigma$. If the condition (\ref{PosMC}) holds and the boundary capacity potential has constant normal derivative on $\Sigma$, then
\begin{eqnarray}\label{DiracCapacityC1}
\lambda_1(\D_\gamma)+\frac{1}{2}\max_{\Sigma}H_g\geq (n-1)\frac{\omega_{n-1}}{|\Sigma|}\mathcal{C}_g(\Sigma,M).
\end{eqnarray}
Moreover, equality holds if, and only if, the Riemannian spin manifold $(M^n,\overline{g})$ with $\overline{g}=\phi^{\frac{4}{n-2}} g$ carries a parallel spinor and $\Sigma$ has positive constant mean curvature. 
\end{corollary} 

Thereafter, we apply Corollary \ref{DiracCapacityConstante} in Section \ref{SSHP} in the context of sub-static manifolds with harmonic potentials which, when combined with the geometric capacitary inequality proved by Agostiniani, Mazzieri and Oronzio \cite{AgostinianiMazzieriOronzio}, leads to a new rigidity result for the Riemannian Schwarzschild manifold (see Corollary \ref{RigidityRS}). Finally, in the last section, we quickly explain how to combine Herzlich's PMT with the approach of Bray and Miao-Hirsch to prove mass-capacity type-inequalities using the spinorial approach (see Theorem \ref{MassCapacityDiracS} and Corollary \ref{MCID}). 


\section{Extrinsic lower bounds for the boundary Dirac operator for singular metrics}\label{ExtrinsicLowerBoundSingular}


In this section, we prove that the estimate (\ref{HMZ}) of Hijazi, Montiel and Zhang which holds for smooth metrics, is also true under less regularity assumptions on the metric. 


\subsection{Singular metrics}\label{SingularMetrics}


Let us now recall the notion of singular metrics following Mantoulidis, Miao and Tam \cite{MantoulidisMiaoTam}. Note that another way to tackle this problem (and even to treat a more general setting) would be to use the spinorial framework developed by Lee and Le Floch \cite{LeeLeFloch} to prove positive mass theorems for manifolds with distributional scalar curvature. 

In the following, a $L^\infty$ Riemannian metric $g$ on a smooth compact $n$-dimensional manifold $\Omega$ with boundary is said to be a {\it singular metric} if $g$ is $C^\infty$ locally away from a compact subset usually called the singular set of $g$ and denoted by $\sing(g)\subset\Omega\setminus\Sigma$. Moreover, the singular set is a disjoint union of compact sets $P$, $Q$ that satisfy:
\begin{enumerate}
\item\label{HyperSing} $P$ is a smoothly embedded two-sided compact hypersurface without boundary,
\begin{enumerate}
\item\label{HyperSinga} near $P$, $g$ can be expressed as 
\begin{eqnarray*}
g(t,z)=dt^2+g_\pm(t,z)
\end{eqnarray*}
for smooth coordinates $(t,z)\in(-a,a)\times P$, $a>0$, where $g_+$ resp. $g_-$ is defined and smooth on $t\geq 0$ resp. $t\leq 0$, and $g_-(0,\cdot)=g_+(0,\cdot)$,
\item\label{HyperSingb} the mean curvature $H_+$, $H_-$ of the unit normal $\frac{\partial}{\partial t}$ at $P$ with respect to $g_+$, $g_-$ satisfy $H_+\leq H_-$. 
\end{enumerate}
\item\label{SobSing} $Q$ is a disjoint union of compact sets $Q_1,\cdots,Q_{N_0}$ such that, for each $k=1,\cdots,N_0$, 
\begin{enumerate}
\item\label{SobSinga} $g$ is $W^{1,q_k}$ in a neighborhood of $Q_k$,
\item\label{codimension} the set $Q_k$ has codimension at least $l_k$, with $l_k>(\frac{2}{n}-\frac{1}{q_k})^{-1}>0$, in the sense that $\limsup_{\varepsilon\rightarrow 0}\varepsilon^{-l_k}|Q_k(\varepsilon)|<\infty$.
\end{enumerate}
\end{enumerate}
Here and below, $c>0$ is a constant independent of $\varepsilon$ which may vary from one line to another and $A(\varepsilon)$ denotes the set of points which are at distance $r<\varepsilon$ from $A$. We are now able to state precisely the main result of this section.
\begin{theorem}\label{HMZ-Singular}
Let $\Omega$ be an $n$-dimensional smooth compact spin manifold with boundary $\Sigma=\coprod_{j=1}^N\Sigma_j$ and let $g$ be a singular metric with singular set $\sing(g)=P\cup Q$ as above. If the scalar curvature $R_g$ of $(\Omega^n,g)$ is nonnegative away from the singular set of $g$ and if $\Sigma$ has positive mean curvature $H_g$ then the first eigenvalue $\lambda_1(\D^{j}_\gamma)$ of the Dirac operator $\D_\gamma$ on each connected component $\Sigma_{j}$ of $\Sigma$ satisfies 
\begin{eqnarray}\label{HMZ-Singular1}
\lambda_{1}(\D^{j}_\gamma)\geq\frac{1}{2}\min_{\Sigma_{j}}H_g
\end{eqnarray} 
for all $j=1,\cdots,N$. Moreover, equality occurs for one $j\in\{1,\cdots,N\}$ if, and only if, $(\Omega^n,g)$ carries a parallel spinor and $\Sigma_j$ has constant mean curvature. If, in addition, $\Omega\setminus Q$ is connected, then $\Sigma$ has to be connected and $H_+=H_-$ along $P$.  
\end{theorem}   
When $\Sigma$ is the boundary of a compact manifold $\Omega$, its mean curvature $H_{\wti{g}}$ will always be computed with respect to the unit normal pointing outside of $\Omega$. 

\begin{remark}
Note that if $Q_k$ is a compact submanifold of codimension at least $2$ for all $k=1,\cdots,N_0$, then $\Omega\setminus Q$ is connected.
\end{remark}

In order to prove Theorem \ref{HMZ-Singular}, we will use the method developed in \cite{MantoulidisMiaoTam} where the authors approximate $g$ by smooth metrics satisfying suitable properties. For simplicity, we may assume that $Q:=Q_1$ with $q:=q_1$ and $l:=l_1$. The properties which are of particular interest for our purpose are listed below. In fact, from \cite[Lemma 3.6]{MantoulidisMiaoTam}, there is $\varepsilon_0>0$ such that if $0<\varepsilon\leq\varepsilon_0$, there is a smooth metric $g_\varepsilon$ on $\Omega$ such that
\begin{enumerate}[label=(\roman*)]
\item $g_{\varepsilon}=g$ outside $P(\varepsilon)\cup Q(\varepsilon)$, $c^{-1}g_\varepsilon\leq g\leq cg_\varepsilon$ for some constant $c>0$ independent of $\varepsilon$,

\item the $W^{1,q}$ norm of $g_\varepsilon$ in $Q(\varepsilon)$ is less that $c$ for some constant $c>0$ independent of $\varepsilon$.
\end{enumerate}
On the other hand, the scalar curvature of a singular metric $g$ is defined away from $\sing(g)$ and so nonnegativity of the scalar curvature only make sense on $\Omega\setminus\sing(g)$. In particular, the Schr\"odinger-Lichnerowicz formula for spinors (see Section \ref{SpinGeometry}), which is at the heart of the proof of the inequality (\ref{HMZ}), does not make sense on the singular support of $g$. Then an important feature of the smooth metric $g_\varepsilon$ is that the scalar curvature $R(g_\varepsilon)$ of $g_\varepsilon$ can still be controlled in term of the scalar curvature of $g$. Indeed, it is proved in \cite[Lemma 3.7]{MantoulidisMiaoTam} that there exists a constant $\alpha>0$ such that if $\varepsilon_0>0$ is small enough so that if $\varepsilon\in(0,\varepsilon_0]$, then for any Lipschitz function $f$ on $\Omega$:
\begin{eqnarray}\label{ScalarCurvatureBoundQ}
\int_{\Omega}R_{g_\varepsilon} f^2d\mu_{g_\varepsilon} & \geq & \int_{\Omega\setminus (P(\varepsilon)\cup Q(\varepsilon))}R_gf^2d\mu_g+c^{-1}\varepsilon^{-2}\int_{P(\frac{1}{400}\varepsilon^2)}(H_--H_+)f^2d\mu_{g_\varepsilon}\nonumber\\& &-c\tau^{\frac{3}{2}}\Big(\int_{Q(\varepsilon)}f^{\frac{2n}{n-2}}d\mu_{g_\varepsilon}\Big)^{\frac{n-2}{n}}-c\tau^{\frac{1}{2}}\int_{\Omega}|\nabla^{g_\varepsilon} f|^2d\mu_{g_\varepsilon}
\end{eqnarray}
where $\tau=\varepsilon^{\alpha/2}$.


\subsection{Generalities on spinors}\label{SpinGeometry}


On an $n$-dimensional manifold $\Omega$ endowed with a spin structure $\Xi$ and a smooth Riemannian metric $g$, there exists a smooth Hermitian vector bundle over $\Omega$ called the spinor bundle and denoted by $\mathbf{S}_g$. The sections of this bundle are called spinors. Moreover, the tangent bundle $T\Omega$ acts on $\mathbf{S}_g$ by Clifford multiplication $X\otimes \psi\mapsto X\cdot\psi$ for any tangent vector fields $X$ and any spinor fields $\psi$. On the other hand, the Riemannian Levi-Civita connection $\nabla^g$ lifts to the so-called spin Levi-Civita connection (also denoted by $\nabla^g$) and defines a metric connection on $\mathbf{S}_g$ that preserves the Clifford multiplication. The Dirac operator is then the first order elliptic differential operator acting on the spinor bundle $\mathbf{S}_g$ locally given by 
\begin{eqnarray*}
D_g:=\sum_{j=1}^n e_j\cdot\nabla^g_{e_j}
\end{eqnarray*}
where $\{e_1,\cdots,e_n\}$ is a $g$-local orthonormal frame of $\Omega$. If $\Omega$ has a boundary $\Sigma:=\partial\Omega$, the spin structure on $\Omega$ induces  a spin structure on its boundary which will be also denoted by $\Xi$. This allows to define the {\it extrinsic} spinor bundle $\RSB_\gamma:=\mathbf{S}_{g|\Sigma}$ over $\Sigma$ on which there exists a metric connection $\nbs^\gamma$. The two spin covariant derivatives are related by the spin Gauss formula which states that 
\begin{eqnarray*}\label{SpinGauss}
\nabla^g_X\varphi=\nbs^\gamma_X\varphi+\frac{1}{2}A_g X\cdot \eta\cdot\varphi
\end{eqnarray*}
for $X\in\Gamma(T\Sigma)$, $\varphi\in\Gamma(\RSB_\gamma)$ and where $A_g:=-\nabla^g\eta$ represents the Weingarten map of $\Sigma$ in $(\Omega^n,g)$ with $\eta$ the unit inner normal to $\Sigma$ in $\Omega$. The extrinsic Dirac operator is defined by taking the Clifford trace of the covariant derivative $\nbs^\gamma$, namely 
\begin{eqnarray*}
\D_\gamma:=\sum_{j=1}^{n-1} e_j\cdot\eta\cdot\nbs^\gamma_{e_j}
\end{eqnarray*}
and is related to the Dirac operator on $\Omega$ by the following formula 
\begin{eqnarray}\label{RelationDirac}
\D_\gamma\varphi=\frac{H_g}{2}\varphi-\eta\cdot D_g\varphi-\nabla_\eta^g\varphi
\end{eqnarray}
which holds for all $\varphi\in\Gamma(\RSB_\gamma)$. From the spin structure $\Xi$ on $\Sigma$, one can also construct an {\it intrinsic} spinor bundle for the induced metric $\gamma$, denoted by $\mathbf{S}_\gamma$, which is also endowed with a spin Levi-Civita connection $\nabla^{\gamma}$. Note that the (intrinsic) Dirac operator $D_{\gamma}$ on $(\Sigma^{n-1},\gamma,\Xi)$ is defined similarly to $D_g$ and $\D_\gamma$. In fact, we have an isomorphism 
$$
\big(\RSB_\gamma,\nbs^\gamma,\D_\gamma\big)\simeq
\left\lbrace
\begin{array}{ll}
\big(\mathbf{S}_{\gamma},\nabla^{\gamma},D_{\gamma}\big) & \text{ if } n \text{ is odd}\\
\big(\mathbf{S}_{\gamma},\nabla^{\gamma},D_{\gamma}\big)\oplus\big(\mathbf{S}_{\gamma},\nabla^{\gamma},-D_{\gamma}\big)& \text{ if } n \text{ is even}
\end{array}
\right.
$$
so that the restriction of a spinor field on $\Omega$ to $\Sigma$ and the extension of a spinor field on $\Sigma$ to $\Omega$ are well-defined. These identifications also imply in particular that the spectrum of the extrinsic Dirac operator is an intrinsic invariant of the boundary: it only depends on the spin and Riemannian structures of $\Sigma$ and not on how it is embedded in $\Omega$. Moreover, the first nonnegative eigenvalue of the extrinsic Dirac operator corresponds to the lowest eigenvalue (in absolute value) of $D_{\gamma}$ and it will be denoted by $\lambda_1(\D_\gamma)$. 

Recall that the Schr\"odinger-Lichnerowicz formula \cite{Lichnerowicz2} which gives a relation between the square of the Dirac operator and the spin Laplacian states that
\begin{eqnarray*}
D_g^2=(\nabla^g\big)^*\nabla^g+\frac{R_g}{4}{\rm Id}_{\mathbf{S}_g}.
\end{eqnarray*}
Once integrated over $\Omega$ (see \cite{HijaziMontielZhang1}), we get 
\begin{eqnarray}\label{Reilly}
\int_{\Omega}\big(|\nabla^g\varphi|^2-|D_g\varphi|^2+\frac{R_g}{4}|\varphi|^2\big)d\mu_g =\int_{\Sigma} \big(\<\D_\gamma\varphi,\varphi\>-\frac{H_g}{2}|\varphi|^2\big)d\mu_\gamma
\end{eqnarray}
for all $\varphi\in\Gamma(\mathbf{S}_g)$ and where $|\,.\,|$ denotes the norm associated to the Hermitian scalar product $\<\,,\,\>$ on $\mathbf{S}_g$.

In the following, the spaces of squared-integrable functions, tensors and spinors on $(\Omega^n,g)$ are denoted by $L^s(g)$ and equipped with the norm
\begin{eqnarray*}
||\varphi||_{L^s(g)}:=\Big(\int_{\Omega}|\varphi|^sd\mu_{g}\Big)^{1/s}
\end{eqnarray*}
for $s\in(0,\infty)$. As well, $W^{1,s}(g)$ denotes the Sobolev spaces of functions, tensors and spinors endowed with the norm
\begin{eqnarray*}\label{SobolevNorms}
||\varphi||_{W^{1,s}(g)}:=||\varphi||_{L^s(g)}+||\nabla^{g}\varphi||_{L^s(g)}.
\end{eqnarray*}

It is a well-known fact that, on a manifold with a fixed spin structure, the spinor bundle depends on the Riemannian metric. We briefly recall here how to identify the spinor bundles over $(\Omega^n,g)$ and $(\Omega^n,g')$ when $g$ and $g'$ are two Riemannian metrics on $\Omega$ using the method of Bourguignon and Gauduchon \cite{BourguignonGauduchon}. In fact, there exists a unique endomorphism $\mathcal{A}^{g}_{g'}$ of $T\Omega$ which is positive, symmetric with respect to $g$ and which maps $g'$-orthonormal frames to $g$-orthonormal frames. It turns out that this map induces an isomorphism $\mathcal{A}^{g}_{g'}:\mathbf{S}_{g'}\rightarrow\mathbf{S}_{g}$ between the spinor bundles which is a fiberwise isometry, compatible with the Clifford multiplication and with inverse $\mathcal{A}^{g'}_{g}$. Then, to compare the corresponding spin Levi-Civita connections, we introduce a third connection 
\begin{eqnarray*}
\wti{\nabla}^{g'}_XY:=\mathcal{A}_{g'}^g\Big(\nabla^{g'}_X\big(\mathcal{A}^{g'}_gY\big)\Big)
\end{eqnarray*}
for all $X$ and $Y$ tangent vectors on $\Omega$ as well as 
\begin{eqnarray*}
\wti{\nabla}^{g'}_X\varphi:=\mathcal{A}_{g'}^g\Big(\nabla^{g'}_X\big(\mathcal{A}^{g'}_g\varphi\big)\Big)
\end{eqnarray*}
for $\varphi\in\Gamma(\BS_g)$, its lift on the spinor bundle $\BS_g$. Now if $(e_i)_{1\leq i\leq n}$ is a local $g$-orthonormal frame on $\Omega$, we get that 
\begin{eqnarray}\label{NablaComp}
\big(\wti{\nabla}_X^{g'}-\nabla_X^g\big)\varphi = \frac{1}{2}\sum_{1\leq k<l\leq n}(\widetilde{\omega}^{g'}_{kl}-\omega^g_{kl})(X)e_k\cdot e_l\cdot\varphi
\end{eqnarray}
where $\widetilde{\omega}^{g'}_{kl}=g(\wti{\nabla}^{g'} e_k,e_l)$ and $\omega^g_{kl}=g(\nabla^g e_k,e_l)$ are the connection $1$-forms associated to $\wti{\nabla}^{g'}$ and $\nabla^g$. In a same way, the Dirac operators $D_g$ and $\wti{D}_{g'}:= \mathcal{A}^g_{g'} \circ D_{g'}\circ\mathcal{A}^{g'}_g$ acting on $\mathbf{S}_g$ can also be related but it is not useful for us here. 


\subsection{Proof of Theorem \ref{HMZ-Singular}}\label{HMZ-SingularMetrics}


Since the metric $g$ is singular, one can consider for all $\varepsilon>0$ the metric $g_\varepsilon$ as in Section \ref{SingularMetrics}. Note that one can find $\wti{\varepsilon}>0$ such that $Q(2\wti{\varepsilon})\cap P(2\wti{\varepsilon})=\emptyset$ and that $Q(2\wti{\varepsilon})\cup P(2\wti{\varepsilon})$ is disjoint from the boundary so that $g_{\varepsilon|\Sigma}$ coincides with $\gamma$ for all $0<\varepsilon\leq\wti{\varepsilon}$. For $(\varepsilon_i)_{i\geq 0}$ a sequence of real numbers with $0<\varepsilon_i\leq \wti{\varepsilon}$ and $\varepsilon_i\rightarrow 0$, we let $g_i:=g_{\varepsilon_i}$. In particular, it holds that $\BS_{g_i|\Sigma}$ coincides with $\RSB_\gamma$ for all $i\geq 0$. On the other hand, it follows from standard arguments (see \cite{HijaziMontielZhang2} for example) that the operator 
\begin{eqnarray*}
D_{g_i}:\big\{\varphi\in\ W^{1,2}(g_i)\,/\, P_\pm\varphi_{|\Sigma}=0\big\}\rightarrow L^2(g_i)
\end{eqnarray*}
is an isomorphism in a such a way that the boundary value problem 
\begin{equation}\label{GeneralBVP}
\left\lbrace
\begin{array}{ll}
D_{g_i}\Phi = 0 & \text { in } \Omega\\
P_\pm\Phi_{|\Sigma}= P_\pm \Psi & \text{ on }\Sigma
\end{array}
\right.
\end{equation}
admits a unique smooth solution $\Phi\in\Gamma(\mathbf{S}_{g_i})$ for every smooth $\Psi\in\Gamma(\RSB_\gamma)$. Here $P_\pm$ denotes the pointwise orthogonal projection defined by
\begin{eqnarray*}
P_\pm:\varphi\in \Gamma(\RSB_{\gamma})\mapsto P_\pm\varphi:=\frac{1}{2}\big(\varphi\pm \sqrt{-1}\eta\cdot\varphi\big)\in \Gamma(V^\pm)
\end{eqnarray*}
where $V^\pm$ is the subbundle of $\RSB_\gamma$ whose fiber is the eigenspace associated with the eigenvalue $\pm 1$ of the involution $\sqrt{-1}\eta\cdot:\RSB_\gamma\rightarrow\RSB_\gamma$. Recall that this maps defines an elliptic boundary condition for $D_{g_i}$ usually referred to as the MIT bag boundary condition.
Now fix $j_0\in\{1,\cdots,N\}$ and let $\phi_i\in\Gamma(\mathbf{S}_{g_i})$ be the unique solution of the boundary value problem (\ref{GeneralBVP}) for $P_+$ with
\begin{eqnarray*}\label{prolonge}
\Psi_{j_0}:=
\left\lbrace
\begin{array}{ll}
\psi_{j_0} & \text{ on } \Sigma_{j_0} \\
0 & \text{ on } \Sigma\setminus\Sigma_{j_0} 
\end{array}
\right.
\end{eqnarray*}
where $\psi_{j_0}$ is an eigenspinor for the Dirac operator on $\Sigma_{j_0}$ associated with the eigenvalue $\lambda_1(\D^{j_0}_\gamma)$. Without loss of generality, we can assume that 
\begin{eqnarray}\label{UnitL2}
\int_{\Omega}|\phi_i|^2d\mu_{g_i}=1
\end{eqnarray}
for all $i\geq 0$. Then the integral version of the Schr\"odinger-Lichnerowicz formula (\ref{Reilly}) on $(\Omega^n,g_i)$ and classical manipulations on the boundary terms (see \cite{Raulot3} for example) yield
\begin{eqnarray}\label{InequalityK0}
\lambda_1(\D_\gamma^{j_0}) \int_{\Sigma_{j_0}}|\phi_i|^2d\mu_\gamma& \geq &  \int_{\Omega}\Big(|\nabla^{g_i}\phi_i|^2+\frac{R_{g_i}}{4}|\phi_i|^2\Big)d\mu_{g_i}+\frac{1}{2}\big(\min_{\Sigma_{j_0}}H_{g}\big)\int_{\Sigma_{j_0}}|\phi_i|^2d\mu_\gamma\nonumber\\
& & +\lambda_1(\D_\gamma^{j_0})\int_{\Sigma_{j_0}}|P_-(\phi_i-\Psi_{j_0})|^2d\mu_\gamma+\int_{\Sigma\setminus\Sigma_{j_0}}H_g|P_-\phi_i|^2d\mu_\gamma.
\end{eqnarray}
Now from (\ref{ScalarCurvatureBoundQ}) with $f=|\phi_i|$, it holds that 
 \begin{eqnarray*}\label{ScalarCurvatureBound}
\int_{\Omega}R_{g_i} |\phi_i|^2d\mu_{g_i} & \geq & \int_{\Omega\setminus (P(\varepsilon_i)\cup Q(\varepsilon_i))}R_g|\phi_i|^2d\mu_g+c^{-1}\varepsilon_i^{-2}\int_{P(\frac{1}{400}\varepsilon^2_i)}(H_--H_+)|\phi_i|^2d\mu_{g_i}\nonumber\\& &-c\tau_i^{\frac{3}{2}}\Big(\int_{Q(\varepsilon_i)}|\phi_i|^{\frac{2n}{n-2}}d\mu_{g_i}\Big)^{\frac{n-2}{n}}-c\tau_i^{\frac{1}{2}}\int_{\Omega}\big|\nabla^{g_i} |\phi_i|\big|^2d\mu_{g_i}.
\end{eqnarray*}
Since the scalar curvature $R_g$ is nonnegative away from $\sing(g)$ and $H_-\geq H_+$ on $P$, we deduce from the previous inequality and the Kato inequality that 
\begin{eqnarray}\label{InequalityK1}
\int_{\Omega}R_{g_i} |\phi_i|^2d\mu_{g_i}\geq-c\tau_i^{\frac{1}{2}}\int_{\Omega}|\nabla^{g_i}\phi_i|^2d\mu_{g_i}-c\tau_i^{\frac{3}{2}}\Big(\int_{Q(\varepsilon_i)}|\phi_i|^{\frac{2n}{n-2}}d\mu_{g_i}\Big)^{\frac{n-2}{n}}.
\end{eqnarray} 
On the other hand, by the Sobolev inequality, the fact that $c g_i\leq g\leq c^{-1}g_i$ and (\ref{UnitL2}), we have 
\begin{eqnarray*}\label{SobolevInequality}
\Big(\int_{\Omega}|\phi_i|^{\frac{2n}{n-2}}d\mu_{g_i}\Big)^{\frac{n-2}{n}}\leq c\Big(\int_{\Omega}|\nabla^{g_i}\phi_i|^2d\mu_{g_i}+1\Big).
\end{eqnarray*}
Combining this fact with (\ref{InequalityK1}), we can rewrite the estimate (\ref{InequalityK0}) as 
\begin{eqnarray}\label{FinalIntegral}
\big(\lambda_1(\D_\gamma^{j_0})-\frac{1}{2}\min_{\Sigma_{j_0}}H_{g}\big)\int_{\Sigma_{j_0}}|\phi_i|^2d\mu_\gamma \geq\big(1-c\tau_i^{\frac{1}{2}}\big)\int_{\Omega}|\nabla^{g_i}\phi_i|^2d\mu_{g_i}-c\tau_i^{\frac{3}{2}}
\end{eqnarray} 
for all $i\geq 0$. Taking the limit as $i$ goes to infinity yields the desired inequality. 

Assume now that equality holds in (\ref{HMZ-Singular1}) and so it follows from (\ref{FinalIntegral}) that 
\begin{eqnarray}\label{ZeroDC}
\lim_{i\rightarrow\infty}\int_{\Omega}|\nabla^{g_i}\phi_i|^2d\mu_{g_i} =0.
\end{eqnarray}
First, remark that our regularity assumptions implies that the metric $g$ must be $C^0\cap W^{1,n}$ and that the sequence $(g_i)$ is uniformly bounded in $W^{1,n}(g)$. Then, it follows from (\ref{NablaComp}) that the ${\rm End}(\mathbf{S}_{g})$-valued one-form
\begin{eqnarray*}
\mathcal{L}^i:X\in\Gamma(TM)\mapsto\mathcal{L}^i_X:=\wti{\nabla}^{g_i}_X-\nabla^{g}_X\in\Gamma\big({\rm End}(\mathbf{S}_{g})\big) 
\end{eqnarray*}
must be $L^n(g)$ and so one compute using the H\"older inequality that  
\begin{eqnarray*}
\int_{\Omega}|\nabla^{g}\varphi|^2d\mu_{g} & \leq &  \int_{\Omega}|\wti{\nabla}^{g_i}\varphi|^2d\mu_g+||\mathcal{L}^i||_{L^n(g)}^2||\varphi||_{L^{n^\ast}(g)}^2
\end{eqnarray*}
that is 
\begin{eqnarray*}
\int_{\Omega}|\nabla^{g}\varphi|^2d\mu_{g} \leq \int_{\Omega}|\wti{\nabla}^{g_i}\varphi|^2d\mu_{g}+c||\varphi||_{L^{n^\ast}(g)}^2
\end{eqnarray*}
for all $\varphi\in\Gamma(\mathbf{S}_g)$ and where $n^\ast=2n/(n-2)$. On the other hand, since the $L^2(g_i)$-norms are uniformly equivalent to the $L^2(g)$-norm, it turns out that the previous inequality implies that 
\begin{eqnarray*}
||\varphi||_{W^{1,2}(g)}\leq c\big(||\mathcal{A}^{g_i}_g\varphi||_{W^{1,2}(g_i)}+||\mathcal{A}^{g_i}_g\varphi||_{L^{n^\ast}(g_i)}\big)
\end{eqnarray*}  
which, from the Sobolev and Kato inequalities, finally yields 
\begin{eqnarray*}
||\varphi||_{W^{1,2}(g)}\leq c||\mathcal{A}^{g_i}_g\varphi||_{W^{1,2}(g_i)}.
\end{eqnarray*}
Combining this estimate with (\ref{UnitL2}) and (\ref{ZeroDC}) show that the sequence $(\Phi_i)$, with $\Phi_i:=\mathcal{A}^g_{g_i}\phi_i\in\Gamma(\mathbf{S}_g)$, is uniformly bounded in $W^{1,2}(g)$. Then, after passing to a subsequence, we conclude that $(\Phi_i)$ converges weakly in $W^{1,2}(g)$, strongly in $L^{2}(g)$ and a.e. on $\Omega$. Denote by $\Phi$ the limit spinor. The strong convergence in $L^{2}(g)$ with (\ref{UnitL2}) and the fact that the sequence $(g_i)$ converges to $g$ in $C^0$ implies that
\begin{eqnarray*}\label{UnitL2Limit}
\int_{\Omega}|\Phi|^2d\mu_{g} = 1.
\end{eqnarray*}
Then the weak convergence in $W^{1,2}(g)$ gives
\begin{eqnarray}\label{Parallel1}
\int_{\Omega}|\nabla^{g}\Phi|^2d\mu_{g}\leq \liminf_{i\rightarrow\infty}\int_{\Omega}|\nabla^{g}\Phi_i|^2d\mu_{g}.
\end{eqnarray} 
Now since $g_i$ and $g$ coincide on $\Omega\setminus\big(Q(\varepsilon_i)\cup P(\varepsilon_i)\big)$, we write
\begin{eqnarray}\label{Parallel2}
\int_{\Omega}|\nabla^g\Phi_i|^2d\mu_g\leq \int_{\Omega}|\wti{\nabla}^{g_i}\Phi_i|^2d\mu_g + \int_{Q(\varepsilon_i)}|\mathcal{L}^{i}\Phi_i|^2d\mu_g+\int_{P(\varepsilon_i)}|\mathcal{L}^{i}\Phi_i|^2d\mu_g.
\end{eqnarray}
The first term in the right-hand side of this inequality tends to zero as $i$ goes to $\infty$ because of (\ref{ZeroDC}). On the other hand, we apply the generalized H\"older inequality in the second term to get  
\begin{eqnarray*}
\int_{Q(\varepsilon_i)}|\mathcal{L}^{i}\Phi_i|^2d\mu_g\leq ||\mathcal{L}^i||_{L^q(g)}^2||\Phi_i||^2_{L^{n^\ast}(g)}\,|Q(\varepsilon_i)|^{2(\frac{1}{n}-\frac{1}{q})}.
\end{eqnarray*}
This term also tends to zero as $i$ goes to infinity since $\mathcal{L}^i$ is uniformly bounded in $L^q(g)$ on $Q(\varepsilon_i)$, $(\Phi_i)$ is uniformly bounded in $L^{n^*}(g)$ by the Sobolev inequality and $Q$ has codimension at least $l$ with $l>nq/(2q-n)$ and $q>n$. Finally, note that on $P(\widetilde{\varepsilon}_0)$, $g$ is equivalent to the Riemannian metric $dt^2+h$ where $h$ is the metric induced by $g$ on $P$ and that $(g_i)$ is uniformly Lipschitz. The last term in the right-hand side of (\ref{Parallel2}) also tends to zero as $i\rightarrow\infty$ since then
\begin{eqnarray*}
\int_{P(\varepsilon_i)}|\mathcal{L}^{i}\Phi_i|^2d\mu_g\leq c||\Phi_i||^2_{L^{n^\ast}(g)}\,\varepsilon_i^{\frac{2}{n}}\leq c\varepsilon_i^{\frac{2}{n}}
\end{eqnarray*}
by the H\"older inequality. We finally have proved that 
\begin{eqnarray*}
\liminf_{i\rightarrow\infty}\int_{\Omega}|\nabla^g\Phi_i|^2d\mu_{g}=0
\end{eqnarray*}
which, with (\ref{Parallel1}), immediately implies that $\Phi$ is a parallel spinor on $(\Omega^n,g)$ as claimed. In particular, since $g$ is smooth away from it singular support, $\Phi$ is smooth on $\Omega\setminus \big(P\cup Q)$. Even more, it follows from \cite[Lemma 3.3]{ShiTam1} that $\Phi$ is H\"older continuous away from $Q$. 

Let us now show that $\Phi_{|\Sigma}=\Psi_{j_0}$. First, the Sobolev trace theorem ensures that the sequence $(\Phi_i)$ converges to $\Phi_{|\Sigma}$ in $L^2(\gamma)$, the space of square integrable spinor fields on $\Sigma$, and so a.e. on $\Sigma$ up to the extraction of a subsequence. Then, since the metric induced on $\Sigma$ by $g_i$ is $\gamma$, it follows that $\Phi_i=\phi_i$ on $\Sigma$ for all $i$. Moreover, from $P_+\phi_{i|\Sigma}=P_+\Psi_{j_0}$ we easily deduce that $P_+\Phi_{|\Sigma}=P_+\Psi_{j_0}$ a.e. on $\Sigma$. On the other hand, we see from (\ref{InequalityK0}) and (\ref{FinalIntegral}) that if equality holds in (\ref{HMZ-Singular1}), we have
\begin{eqnarray*}
\lim_{i\rightarrow\infty}\int_\Sigma |P_-\Phi_i-P_-\Psi_{j_0}|^2d\mu_\gamma=0.
\end{eqnarray*}
This means that $P_-\phi_{i|\Sigma}$ converges to $P_-\Psi_{j_0}$ in $L^2(\gamma)$ and so a.e. on $\Sigma$ up to a subsequence. In particular, $P_-\Phi_{|\Sigma}=P_-\Psi_{j_0}$ a.e. on $\Sigma$ which implies that $\Phi_{|\Sigma}=\Psi_{j_0}$ a.e. on $\Sigma$. However, since both of these spinor fields are continuous on $\Sigma$, this equality holds everywhere on $\Sigma$. Now from the formula (\ref{RelationDirac}), we get
\begin{eqnarray*}
\frac{1}{2}H_g\Phi_{|\Sigma_{j_0}}=\D_\gamma\Phi_{|\Sigma_{j_0}}=\D_\gamma\psi_{j_0}=\lambda_{1}(\D^{j_0}_\gamma)\psi_{j_0}=\lambda_{1}(\D^{j_0}_\gamma)\Phi_{|\Sigma_{j_0}}
\end{eqnarray*}
and so $H_g$ is constant on $\Sigma_{j_0}$. 

Now if $\Omega\setminus Q$ is connected then, since $\Phi$ is parallel on $\Omega$ and continuous on $\Omega\setminus Q$, it has a positive constant norm away from $Q$. This is impossible if $\Sigma$ has more than one component since otherwise it should be zero on $\Sigma\setminus\Sigma_{j_0}$. It remains to prove that $H_-=H_+$ on $P$. For this, we note that for all spinor fields $\varphi$ which are smooth on $P(\wti{\varepsilon})\setminus P$ and $W^{1,2}(g)$ on $P(\wti{\varepsilon})$, we have
\begin{eqnarray*}
\int_{P(\widetilde{\varepsilon})}|D_g\varphi|^2 & = & \int_{P(\widetilde{\varepsilon})}|\nabla^g\varphi|^2+\int_P(H_--H_+)|\varphi|^2\\
& + & \int_{P_{\wti{\varepsilon}}}\< D_g\varphi-\frac{\partial}{\partial t}\cdot\nabla^g\varphi,\varphi\>+\int_{P_{-\wti{\varepsilon}}}\<D_g\varphi+\frac{\partial}{\partial t}\cdot\nabla^g\varphi,\varphi\>
\end{eqnarray*}
where $P_t$ denotes the hypersurface at oriented distance $t$ of $P$ in the tubular neighborhood $P(\wti{\varepsilon})$. Applying this formula to the parallel spinor $\Phi$ leads to the fact that
\begin{eqnarray*}
\int_P(H_--H_+)=0
\end{eqnarray*}
and so $H_-=H_+$ since $H_-\geq H_+$ on $P$.

Conversely, assume that $\Omega$ endowed with a singular metric $g$ carries a parallel spinor fields $\Phi\in\Gamma(\mathbf{S}_{g})$ and has a boundary $\Sigma$ with positive constant mean curvature $H_{g}$. Then, it is straightforward to compute using (\ref{RelationDirac}) that the restriction of $\Phi$ to $\Sigma$ is an eigenspinor for the Dirac operator $\D_\gamma$ associated with the eigenvalue $H_{g}/2$. On the other hand, since $(\Omega^n,g)$ has a parallel spinor it has to be Ricci-flat away from $P\cup Q$ and so the inequality (\ref{HMZ-Singular1}) applies. This leads to the conclusion that it is actually an equality as desired.


\section{Proof of Theorem \ref{DiracEigenvalueBound} and its corollaries}


{\it Proof of Theorem \ref{DiracEigenvalueBound}.} Consider the metric conformally related to $g$ defined by $\overline{g}=\phi^{\frac{4}{n-2}}g$ on $M$ where $\phi$ is the boundary capacity potential satisfying (\ref{Capacity}). From the classical relation between scalar curvatures in the same conformal class, it holds that
\begin{eqnarray*}
R_{\overline{g}}=\phi^{-\frac{n+2}{n-2}}\Big(-4\frac{n-1}{n-2}\Delta_g \phi+R_g\phi\Big)\geq 0
\end{eqnarray*}
since $\phi$ is harmonic and $R_g$ is nonnegative. Now it follows from \cite{MantoulidisMiaoTam,MiaoHirsch} that we can invert the coordinates at infinity with the help of the Kelvin transform to compactify the manifold $M$ by adding a point $p_\infty$ at infinity to a get a smooth manifold $M_\infty$ for which the metric $\overline{g}$ extends to a $W^{1,q}$ metric for some $q>n$ at $p_\infty$ . Denote by $\overline{g}_\infty$ the extended metric and remark that this metric is smooth on $M$ because $\overline{g}_\infty=\overline{g}$ on $M$. Moreover, since the manifold $M$ is spin, its one point compactification $M_\infty$ is also spin and the spin structure induced on $\Sigma$ remains unchanged. On the other hand, the mean curvature of $\Sigma$ for the metric $\overline{g}$ can be easily seen to be 
\begin{eqnarray}\label{ConfMC}
H_{\overline{g}}=\phi^{-\frac{2}{n-2}}\Big(-2\frac{n-1}{n-2}\frac{\partial \phi}{\partial\nu}-H_g\phi\Big)=-2\frac{n-1}{n-2}\frac{\partial \phi}{\partial\nu}-H_g>0
\end{eqnarray}
because of (\ref{PosMC}) and since $\phi=1$ on $\Sigma$. For the same reason, the metric $\overline{g}$ restricted to $\Sigma$ coincides with $\gamma$. The metric $\overline{g}_\infty$ being singular in the sense of Section \ref{SingularMetrics} with $P=\emptyset$ and $Q=\{p_\infty\}$, we can apply Theorem \ref{HMZ-Singular} and then the inequality (\ref{DiracCapacity}) follows directly from (\ref{HMZ-Singular1}) and (\ref{ConfMC}). Since $M_\infty\setminus\{p_\infty\}=M$ is connected, the equality case is also a direct consequence of the equality case of Theorem \ref{HMZ-Singular}. 
\hfill$\square$

\vspace{0.2cm}

{\it Proof of Corollary \ref{EqualityCaseThree}.} Assume here that $(M^3,g)$ is an $3$-dimensional, complete, asymptotically flat manifold with nonnegative scalar curvature whose inner boundary is a union of $2$-spheres for which equality holds in (\ref{DiracCapacity}). It follows from the equality case in Theorem \ref{DiracEigenvalueBound} that $(M^3,\ov{g})$ carries a parallel spinor and that $\Sigma$ is a $2$-sphere with positive constant mean curvature. In particular, the manifold $(M^3,\ov{g})$ is Ricci flat and so it is flat as $M$ is $3$-dimensional. In fact, it follows from \cite{SmithYang} that the metric $\ov{g}_\infty$ defined on $M^3_\infty$ in the previous proof, extends smoothly (after perhaps a change of smooth structure) across $p_\infty$ as the singular support of $\ov{g}_\infty$ is reduced to an isolated point. Then we deduce that $(M^3_\infty,\ov{g}_\infty)$ is a handlebody with a flat metric and a connected mean convex spherical boundary, that is a $3$-dimensional ball. The end of the proof proceeds then exactly as in \cite[p. 26]{MantoulidisMiaoTam} to conclude that $(M^3,g)$ is isometric to the exterior of a coordinate sphere in a Riemannian Schwarzschild manifold.
\hfill$\square$

\vspace{0.2cm}

{\it Proof of Corollary \ref{MinimalEstimate}.} The inequality (\ref{DiracCapacity}) applies with $H_g=0$ and so (\ref{DiracCapacityMinimal}) is true. The equality case follows directly from Theorem \ref{DiracEigenvalueBound}. 
\hfill$\square$

\vspace{0.2cm}

{\it Proof of Corollary \ref{DiracCapacityConstante}.} By assumption, the normal derivative of the boundary capacity potential $\phi$ of $\Sigma$ in $(M^n,g)$ is constant and, since $\Sigma$ is connected, the second equality in (\ref{DefCap}) shows that 
\begin{eqnarray*}
\mathcal{C}_g(\Sigma,M)=-\frac{|\Sigma|}{(n-2)\omega_{n-1}}\frac{\partial\phi}{\partial \nu}.
\end{eqnarray*}
The inequality (\ref{DiracCapacityC1}) as well as the corresponding equality case follows from Theorem \ref{DiracEigenvalueBound}. 
\hfill$\square$


\section{Sub-static manifolds with harmonic potentials}\label{SSHP}


In \cite[Theorem 1.1]{AgostinianiMazzieriOronzio}, Agostiniani, Mazzieri and Oronzio proved a nice inequality similar to the Penrose inequality in the context of sub-static manifolds with harmonic potentials. A triple $(M^n,g,u)$ is said to be sub-static harmonic if $(M^n,g)$ is a smooth, connected, complete asymptotically flat, $n$-dimensional Riemannian manifold, with $n\geq 3$ and with smooth compact boundary $\Sigma$ and $u\in C^\infty(M)$ satisfies the system
\begin{equation}\label{SubStaticTriple}
\left\lbrace
\begin{array}{ll}
u Ric_g-\nabla^g du\geq 0 & \text{ in } M\\
\Delta_g u=0 & \text{ in } M\\ 
u=0 & \text{ on } \Sigma \\
u\rightarrow 1 & \text{ as } |x|\rightarrow\infty.
\end{array}   
\right.
\end{equation}
Here $Ric_g$ is the Ricci curvature tensor of $(M^n,g)$. In other words, a triple $(M^n,g,u)$ is a sub-static triple if, and only if, the boundary capacity potential $\phi=1-u$  of $\Sigma$ in $(M^n,g)$ satisfies $(1-\phi) Ric_g+\nabla^g d\phi\geq 0$. In particular, it follows from (\ref{ExpansionBCP}) that 
\begin{eqnarray*}
u(x)=1-\frac{\mathcal{C}_g(\Sigma,M)}{r^{n-2}}+o_2(r^{2-n})
\end{eqnarray*}
as $r\rightarrow\infty$. In this context, they proved:
\begin{theorem}$($\cite{AgostinianiMazzieriOronzio}$)$\label{CapacityPenrose}
Let $(M^n,g,u)$ be a sub-static harmonic triple with associated capacity $\mathcal{C}_g(\Sigma,M)$ and suppose that $\Sigma$ is connected. Then
\begin{eqnarray}\label{CapacityPenrose1}
\mathcal{C}_g(\Sigma,M)\geq\frac{1}{2}\Big(\frac{|\Sigma|}{\omega_{n-1}}\Big)^{\frac{n-2}{n-1}}.
\end{eqnarray}
Moreover, the equality holds if, and only if, $(M^n,g)$ is isometric to $\big(\mathbb{M}^n_m(r_m),g_m\big)$ with $r_m=(m/2)^{1/(n-2)}$ and $m=\mathcal{C}_g(\Sigma,M)$.
\end{theorem}
From the two first conditions in (\ref{SubStaticTriple}) it holds that $\nabla^g du=0$ on $\Sigma$. Then, since $\Sigma$ is connected, we get from the 
Hopf lemma that $|\nabla^g u|$ is a positive constant on $\Sigma$ and therefore $|\nabla^g\phi|$ too. In particular, Corollary \ref{DiracCapacityConstante} applies in this situation and leads to the following estimate when combined with Theorem \ref{CapacityPenrose}.
\begin{corollary}\label{DiracPenroseInequality}
Let $(M^n,g,u)$ be a spin sub-static harmonic triple and suppose that $\Sigma$ is connected. Then
\begin{eqnarray}\label{DiracPenrose}
\lambda_1(\D_\gamma)\geq\frac{n-1}{2}\Big(\frac{\omega_{n-1}}{|\Sigma|}\Big)^{\frac{1}{n-1}}.
\end{eqnarray}
Moreover, the equality holds if, and only if, $(M^n,g)$ is isometric to $\big(\mathbb{M}^n_m(r_m),g_m\big)$ with $r_m=(m/2)^{1/(n-2)}$ and 
\begin{eqnarray*}
m=\frac{\lambda_1(\D_\gamma)}{n-1}\frac{|\Sigma|}{\omega_{n-1}}.
\end{eqnarray*}
\end{corollary}
For $n=3$, this inequality is nothing else but the B\"ar-Hijazi inequality (\ref{BarHijazi}) and this is in fact not surprising because of the black hole uniqueness theorem for sub-static harmonic triple \cite[Theorem 1.2]{AgostinianiMazzieriOronzio}. Indeed, this result asserts that any such triple for which there exists a chart at infinity with $R_g=O(r^{-q})$ for some $q>n$ has to be isometric to $\big(\mathbb{M}^n_m(r_m),g_m\big)$ with $r_m=(m/2)^{1/(n-2)}$. In particular, $\Sigma$ is a round sphere and then the inequality (\ref{DiracPenrose}) is in fact an equality (and so is (\ref{CapacityPenrose1})). It remains an open question to see whether it is possible to remove this assumption. For $n\geq 4$, if stronger conditions on the asymptotic behaviors of the metric $g$ and the harmonic potential are made and if $M$ is spin, similar rigidity results can be deduced from \cite{Raulot12}. In the situation here, one can obtain from Corollary \ref{DiracPenroseInequality} the following uniqueness result:
\begin{corollary}\label{RigidityRS}
Let $(M^n,g,u)$ be a spin sub-static harmonic triple whose boundary is isometric to a round sphere with radius $R>0$. Then $(M^n,g)$ is isometric to $\big(\mathbb{M}^n_m(r_m),g_m\big)$ with $r_m=(m/2)^{1/(n-2)}$ and $m=R^{n-2}/2$.
\end{corollary}
This follows directly from the fact that the first eigenvalue of the Dirac operator on a round sphere with radius $R$ is $\lambda_1(\D_\gamma)=(n-1)/(2R)$ so that the equality occurs in (\ref{DiracPenrose}) under the assumptions of Corollary \ref{RigidityRS}. Note that, unlike the above results, we made not use of the PMT in our proof.

\section{Mass-capacity inequalities}\label{PMT-H} 


In this last part, we briefly explain that the mass-capacity inequality (\ref{MCI}) is true if we assume that some estimates regarding $\lambda_1(\D_\gamma)$ holds on $\Sigma$. More precisely, we prove:
\begin{theorem}\label{MassCapacityDiracS}
Let $(M^n,g)$ be an $n$-dimensional, complete, spin asymptotically flat manifold with nonnegative scalar curvature and with a compact inner boundary $\Sigma$ such that 
\begin{eqnarray}\label{HerzlichCapacity1}
\lambda_1(\D_\gamma)\geq\frac{1}{2}\max_{\Sigma}\Big(\frac{2\alpha}{\alpha-2}\frac{n-1}{n-2}\frac{\partial\phi}{\partial\nu}+H_g\Big)
\end{eqnarray}
for some $\alpha\in[0,2[$. Then the ADM mass $m_{ADM}(M,g)$ of $(M^n,g)$ satisfies $$m_{ADM}(M,g)\geq\alpha\,\mathcal{C}_g(\Sigma,M).$$
\end{theorem}

{\it Proof:} Consider the metric $g_\alpha=\phi_\alpha^{\frac{4}{n-2}} g$ conformally related to $g$ and where $\phi_\alpha$ is the smooth positive function defined on $M$ by
$$
\phi_\alpha:=1-\frac{\alpha}{2}\phi>0
$$ 
for $\alpha\in[0,2[$. From the asymptotic expansion (\ref{ExpansionBCP}) of $\phi$, it is immediate to check that $(M^n,g_\alpha)$ is an asymptotically flat manifold with ADM mass given by
\begin{eqnarray}\label{MassConformal}
m_{ADM}(M,g_\alpha)=m_{ADM}(M,g)-\alpha\,\mathcal{C}_g(\Sigma,M).
\end{eqnarray}
Moreover, the scalar curvature of $(M^n,g_\alpha)$ is easily computed to be
\begin{eqnarray*}
R_{g_\alpha}=\phi^{-\frac{n+2}{n-2}}_\alpha\Big(-4\frac{n-1}{n-2}\Delta_g\phi_\alpha+R_g\phi_\alpha\Big)\geq 0
\end{eqnarray*}
since $\phi_\alpha$ is harmonic and $R_g$ is nonnegative. As well, the mean curvature of $\Sigma$ in $(M^n,g_\alpha)$ is 
\begin{eqnarray*}
H_{g_\alpha}=\phi_\alpha^{-\frac{n}{n-2}}\Big(2\frac{n-1}{n-2}\frac{\partial \phi_\alpha}{\partial\nu}+H_g\phi_\alpha\Big)
\end{eqnarray*}
which can be rewritten as
\begin{eqnarray}\label{MeanCurvatureAlpha}
H_{g_\alpha}=\Big(\frac{2}{2-\alpha}\Big)^{\frac{2}{n-2}}\Big(\frac{2\alpha}{\alpha-2}\frac{n-1}{n-2}\frac{\partial\phi}{\partial\nu}+H_g\Big).
\end{eqnarray}
Now, observe that since the metrics $g_\alpha$ and $g$ are homothetic when restricted to $\Sigma$, it is easy to check that their corresponding first Dirac eigenvalues satisfy
\begin{eqnarray*}
\lambda_1(\D_{\gamma_\alpha})=\Big(\frac{2}{2-\alpha}\Big)^{\frac{2}{n-2}}\lambda_1(\D_\gamma)
\end{eqnarray*}
where $\lambda_1(\D_{\gamma_\alpha})$ corresponds to the first eigenvalue of the Dirac operator on $\Sigma$ endowed with $\gamma_\alpha$, the metric induced by $g_\alpha$. This fact together with (\ref{HerzlichCapacity1}) and (\ref{MeanCurvatureAlpha}) allows to conclude that 
\begin{eqnarray*}
\lambda_1(\D_{\gamma_\alpha})\geq\frac{1}{2}\max_{\Sigma} H_{g_\alpha}.
\end{eqnarray*}
We thus have shown that the manifold $(M^n,g_\alpha)$ satisfies the assumptions of Theorem \ref{HerzlichPMT-n} and then its ADM mass $m_{ADM}(M,g_\alpha)$ is nonnegative. The conclusion follows directly from (\ref{MassConformal}).  
\hfill$\square$

\begin{remark}
It follows from the estimate (\ref{DiracCapacity}) that the previous mass-capacity inequality holds if 
\begin{eqnarray*}
\min_{\Sigma}\Big(-2\frac{n-1}{n-2}\frac{\partial\phi}{\partial\nu}-H_g\Big) \geq \max_\Sigma\Big(\frac{2\alpha}{\alpha-2}\frac{n-1}{n-2}\frac{\partial\phi}{\partial\nu}+H_g\Big).
\end{eqnarray*} 
However, it is straightforward to see that this assumption
 implies the pointwise condition (\ref{CondMH}). This mean that, in this situation, Theorem \ref{MassCapacityDiracS} is a direct consequence of the works of Bray \cite{Bray} and Hirsch and Miao \cite{MiaoHirsch}. 
\end{remark}

As a direct application of Theorem \ref{MassCapacityDiracS}, we obtain a natural generalization of a PMT for manifolds with boundary \cite{Herzlich1,Herzlich2}. Indeed, combining this result with the B\"ar-Hijazi inequality 
\begin{eqnarray}\label{BarHijazi}
\lambda_1(\D_\gamma)\geq 2\sqrt{\frac{\pi}{|\Sigma|}},
\end{eqnarray}
for $n=3$ and the Hijazi inequality  
\begin{eqnarray}\label{HijaziInequality}
\lambda_1(\D_\gamma)\geq\frac{1}{2}|\Sigma|^{-\frac{1}{n-1}}\sqrt{\frac{n-1}{n-2}\mathcal{Y}\big(\Sigma,[\gamma]\big)}
\end{eqnarray}
for $n\geq 4$ and where $\mathcal{Y}\big(\Sigma,[\gamma]\big)$ is the Yamabe invariant of $(\Sigma^{n-1},\gamma)$, yields the following version of the mass-capacity inequality:
\begin{corollary}\label{MCID}
Let $(M^n,g)$ be an $n$-dimensional, complete, spin asymptotically flat manifold with nonnegative scalar curvature and with a compact inner boundary $\Sigma$. Assume that one of the following condition holds on $\Sigma$ for some $\alpha\in[0,2[$:  
\begin{enumerate}[label=(\roman*)]
\item $n=3$, $\Sigma$ is a topological sphere and 
\begin{eqnarray*}
H_g-\frac{4\alpha}{2-\alpha}\frac{\partial\phi}{\partial\nu}\leq 4\sqrt{\frac{\pi}{|\Sigma|}}
\end{eqnarray*}
\item $n\geq 4$, $\Sigma$ is a manifold with positive Yamabe invariant and 
\begin{eqnarray*}
H_g-\frac{2\alpha}{2-\alpha}\frac{n-1}{n-2}\frac{\partial\phi}{\partial\nu}\leq |\Sigma|^{-\frac{1}{n-1}}\sqrt{\frac{n-1}{n-2}\mathcal{Y}\big(\Sigma,[\gamma]\big)}.
\end{eqnarray*}
\end{enumerate}
Then the ADM mass $m_{ADM}(M,g)$ of $(M^n,g)$ satisfies $$m_{ADM}(M,g)\geq\alpha\,\mathcal{C}_g(\Sigma,M).$$
Moreover, equality occurs if, and only if, $(M^n,g)$ is isometric to the exterior region outside a rotationally symmetric sphere in the Riemannian Schwarzschild manifold of mass $\alpha\,\mathcal{C}_g(\Sigma,M)$. 
\end{corollary}

{\it Proof:} The inequality is a direct consequence of a combination of Theorem \ref{MassCapacityDiracS} with (\ref{BarHijazi}) and (\ref{HijaziInequality}). Assume that equality holds. If $\alpha=0$, it is nothing else that the PMT of Herzlich. If $\alpha>0$ and equality occurs, it is easy to observe from the proof of Theorem \ref{MassCapacityDiracS} that $(M^n,g_\alpha)$ satisfies the equality case of the PMT \cite[Proposition 2.1]{Herzlich1} and \cite[Proposition 2.1]{Herzlich2} and so it has to be isometric to the exterior of a round sphere with radius $r_\alpha>0$ in the Euclidean space. Then it follows that the smooth function $\phi_\alpha^{-1}$ satisfies
$$
\left\lbrace
\begin{array}{ll}
\Delta_\delta\phi_\alpha^{-1}= 0 & \text{ in } \mathbb{R}^n\setminus B(0,r_\alpha) \\
\phi_{\alpha}^{-1}=2/(2-\alpha) & \text{ on } \mathbb{S}^{n-1}_{r_\alpha}\\ 
\phi_{\alpha}^{-1}\rightarrow 1 & \text{ as } |x|\rightarrow\infty
\end{array}
\right.
$$
where $\Delta_\delta$ is the Laplace operator in the Euclidean space. It turns out that the unique solution of the aforementioned boundary problem is easily seen to be 
\begin{eqnarray*}
\phi_\alpha^{-1}(x)=1+\frac{m}{2r^{n-2}}
\end{eqnarray*}
with $m=2\alpha r_\alpha^{n-2}/(2-\alpha)$ and thus $(M^n,g)$ is isometric to $\big(\mathbb{M}^n_m(r_\alpha),g_m\big)$. The converse statement follows from the computations of Appendix \ref{RiemSchwar}. 
\hfill$\square$

\begin{remark}
If the boundary is not connected, the condition (i) or (ii) in Corollary \ref{MCID} as well as the condition (\ref{HerzlichCapacity1}) in Theorem \ref{MassCapacityDiracS} is assumed to hold on each connected components of $\Sigma$.  
\end{remark}


\appendix



\section{Mass and capacity}\label{SectionCapacity}


A smooth, connected, $n$-dimensional Riemannian manifold $(M^n,g)$ is said to be asymptotically flat (of order $p$) if there exists a compact subset $K\subset M$ such that $M\setminus K$ is a finite disjoint union of ends $M_k$, each of them being diffeomorphic to $\mathbb{R}^n$ minus a closed ball $\overline{B}$ by a coordinate chart in which the components of the metric satisfy 
\begin{eqnarray*}
 g_{ij} =  \delta_{ij} +O_2\big(r^{-p}\big)\quad\text{and}\quad R_g\in L^1(M)
\end{eqnarray*}
for $i,j=1,...,n$ and $p>(n-2)/2$. Here $O_2(r^{-p})$ refers to a real-valued function $f$ such that
\begin{eqnarray*}
|f(x)|+r|\partial f(x)|+r^2|\partial^2f(x)|\leq Cr^{-p}
\end{eqnarray*}
as $r$ goes to infinity, for some constant $C>0$ and where $\partial$ is the standard derivative in the Euclidean space. Such a coordinate chart is often referred to as a chart at infinity. In the following, we assume that $k=1$ and the general case can be treated in a similar way.  

On an asymptotically flat manifold $(M^n,g)$, the ADM mass is defined by 
\begin{eqnarray*}
m_{ADM}(M,g):=\frac{1}{2(n-1) \omega_{n-1}}\lim_{R\rightarrow+\infty}\sum_{i,j=1}^n\int_{\mathcal{S}_R}(g_{ij,i}-g_{ii,j}) \frac{x^j}{r}d\overline{\mu}_{{\mathcal{S}_R}}
\end{eqnarray*}
where $\mathcal{S}_R$ stands for a coordinate sphere of radius $R>0$, $d\overline{\mu}_{{\mathcal{S}_r}}$ its Euclidean Riemannian volume element and $g_{ij,s}$ the derivative of the metric components in the coordinate chart. This definition seems to depend on a particular choice of the coordinates chart, however, as proved independently by Bartnik \cite{Bartnik1} and Chru\'sciel \cite{Chrusciel1}, it is a well-defined geometric invariant. 

On asymptotically flat manifolds of order $p$, the boundary capacity potential $\phi\in C^\infty(M)$ satisfying (\ref{Capacity}) has the following expansion 
\begin{eqnarray}\label{ExpansionBCP}
\phi(x)=\frac{\mathcal{C}_g(\Sigma,M)}{r^{n-2}}+O_2(r^{-(n-2+p)})
\end{eqnarray}
as $r\rightarrow\infty$ (see \cite[Theorem 2.2]{AgostinianiMazzieriOronzio} for example).


\section{The Riemannian Schwarzschild manifolds}\label{RiemSchwar}


Here we recall some standard computations in the Riemannian Schwarzschild manifolds. The Riemannian Schwarzschild manifold of mass $m\in\mathbb{R}$ is the Riemannian manifold $(\mathbb{M}^n_m,g_m)$ where 
$$
\mathbb{M}^n_m:=
\left\lbrace
\begin{array}{ll}
\mathbb{R}^n\setminus\{0\} & \text{ if } m>0\\
\mathbb{R}^n & \text{ if } m=0\\
\mathbb{R}^{n}\setminus\big\{r\leq\big(|m|/2\big)^{1/(n-2)}\big\} & \text{ if } m<0
\end{array}
\right.
$$
and with metric 
\begin{eqnarray*}\label{SchwarzschildMetric}
g_m=\Big(1+\frac{m}{2r^{n-2}}\Big)^{\frac{4}{n-2}}\delta
\end{eqnarray*}
where $r:=|x|$ is the Euclidean radius for $x\in\mathbb{M}^n_M$. It is a static manifold in the sense that the Lorentzian manifold 
\begin{eqnarray*}
\Big(\mathfrak{L}^{n+1}:=\mathbb{R}\times\mathbb{M}^n_m,\mathfrak{g}_m:=-N_m^2\,dt^2+g_m\Big)
\end{eqnarray*}
is a spacetime which satisfies the Einstein vacuum equations with zero cosmological constant. Here $N_m$ denotes the smooth harmonic function given by
\begin{eqnarray*}
N_m(x)=\big(1-\frac{m}{2r^{n-2}}\big)\big(1+\frac{m}{2r^{n-2}}\big)^{-1}
\end{eqnarray*}
generally referred to as the lapse function. For $r_0\in(r_\ast,\infty)$ with $r_\ast=0$ if $m\geq 0$ and $r_\ast=(|m|/2)^{1/(n-2)}$, we consider the exterior of the region outside a rotationally symmetric sphere defined by
\begin{eqnarray*}
\mathbb{M}_m^n(r_0):=\Big\{x\in\mathbb{M}^n_m\,/\,r\geq r_0\Big\}.
\end{eqnarray*}
This is an $n$-dimensional, spin, complete, asymptotically flat manifold with zero scalar curvature and connected inner boundary $\Sigma_{r_0}$ with induced metric $\gamma_{r_0}$ isometric to a round sphere with radius
\begin{eqnarray*}\label{RadiusSS}
r_{g_m,r_0}:=r_0\big(1+\frac{m}{2r_0^{n-2}}\big)^{\frac{2}{n-2}}
\end{eqnarray*}
and constant mean curvature 
\begin{eqnarray*}\label{MeanCurvatureSS}
H_{g_m,r_0}=\frac{n-1}{r_0}\big(1-\frac{m}{2r_0^{n-2}}\big)\big(1+\frac{m}{2r_0^{n-2}}\big)^{-\frac{n}{n-2}}.
\end{eqnarray*}
Remark that for $m>0$, the region $\mathbb{R}\times\mathbb{M}^n_m(r_m)$ with $r_m=(m/2)^{1/(n-2)}$ represents the exterior of a black hole with event horizon at $r=r_m$ in $(\mathfrak{L}^{n+1},\mathfrak{g}_m)$. On the other hand, the boundary capacity potential of $\Sigma_{r_0}$ in $\big(\mathbb{M}^n_m,g_m\big)$ can be computed to be 
\begin{eqnarray*}
\phi_{r_0}(x)=\big(1+\frac{m}{2r_0^{n-2}}\big)\big(1+\frac{m}{2r^{n-2}}\big)^{-1}\Big(\frac{r_0}{r}\Big)^{n-2}.
\end{eqnarray*}
Thus it holds on $\Sigma_{r_0}$ that 
\begin{eqnarray*}\label{NormalDerivativeBCP}
\frac{\partial\phi_{r_0}}{\partial\nu}=-\frac{n-2}{r_0}\big(1+\frac{m}{2r_0^{n-2}}\big)^{-\frac{n}{n-2}}
\end{eqnarray*}
and then
\begin{eqnarray*}\label{SchwarzschildCapacity}
\mathcal{C}_{g_m}\big(\Sigma_{r_0},\mathbb{M}^n_m(r_0)\big)=\frac{m}{2}+r_0^{n-2}.
\end{eqnarray*}
It is also relevant to note that 
\begin{eqnarray*}
-2\frac{n-1}{n-2}\frac{\partial\phi_{r_0}}{\partial\nu}-H_{g_m,r_0}=\frac{n-1}{r_{g_m,r_0}}
\end{eqnarray*}
so that this leads to the following interpretation 
\begin{eqnarray*}
\lambda_1(\D_{r_0})=-\frac{n-1}{n-2}\frac{\partial\phi_{r_0}}{\partial\nu}-\frac{H_{g_m,r_0}}{2}
\end{eqnarray*}
in terms of $\lambda_1(\D_{r_0})$ the first eigenvalue of the Dirac operator $\D_{r_0}$ of $(\Sigma_{r_0},\gamma_{r_0})$. This implies that the manifold $\mathbb{M}_m^n(r_0)$ satisfies the equality case in the estimate of Theorem \ref{DiracEigenvalueBound}.


\bibliographystyle{alpha}     
\bibliography{BiblioHabilitation}

\begin{thebibliography}{MMT20}

\bibitem[AMO22]{AgostinianiMazzieriOronzio}
V.~Agostiniani, L.~Mazzieri, and F.~Oronzio.
\newblock A geometric capacitary inequality for sub-static manifolds with
  harmonic potentials.
\newblock {\em Math. Eng.}, 4(2):Paper No. 013, 40, 2022.

\bibitem[Bar86]{Bartnik1}
R.~Bartnik.
\newblock The mass of an asymptotically flat manifold.
\newblock {\em Comm. Pure Appl. Math.}, 39(5):661--693, 1986.

\bibitem[BG92]{BourguignonGauduchon}
J.-P. Bourguignon and P.~Gauduchon.
\newblock Spineurs, op\'{e}rateurs de {D}irac et variations de m\'{e}triques.
\newblock {\em Comm. Math. Phys.}, 144(3):581--599, 1992.

\bibitem[BL09]{BrayLee}
H.~L. Bray and D.~A. Lee.
\newblock On the {R}iemannian {P}enrose inequality in dimensions less than
  eight.
\newblock {\em Duke Math. J.}, 148(1):81--106, 2009.

\bibitem[Bra01]{Bray}
H.~L. Bray.
\newblock Proof of the {R}iemannian {P}enrose inequality using the positive
  mass theorem.
\newblock {\em J. Differential Geom.}, 59(2):177--267, 2001.

\bibitem[Chr86]{Chrusciel1}
P.~T. Chru\'{s}ciel.
\newblock Boundary conditions at spatial infinity from a {H}amiltonian point of
  view.
\newblock In {\em Topological properties and global structure of space-time
  ({E}rice, 1985)}, volume 138 of {\em NATO Adv. Sci. Inst. Ser. B Phys.},
  pages 49--59. Plenum, New York, 1986.

\bibitem[Her97]{Herzlich1}
M.~Herzlich.
\newblock A {P}enrose-like inequality for the mass of {R}iemannian
  asymptotically flat manifolds.
\newblock {\em Comm. Math. Phys.}, 188(1):121--133, 1997.

\bibitem[Her02]{Herzlich2}
M.~Herzlich.
\newblock Minimal surfaces, the {D}irac operator and the {P}enrose inequality.
\newblock In {\em S\'{e}minaire de {T}h\'{e}orie {S}pectrale et
  {G}\'{e}om\'{e}trie, {V}ol. 20, {A}nn\'{e}e 2001--2002}, volume~20 of {\em
  S\'{e}min. Th\'{e}or. Spectr. G\'{e}om.}, pages 9--16. Univ. Grenoble I,
  Saint-Martin-d'H\`eres, 2002.

\bibitem[HI01]{HuiskenIlmanen}
G.~Huisken and T.~Ilmanen.
\newblock The inverse mean curvature flow and the {R}iemannian {P}enrose
  inequality.
\newblock {\em J. Differential Geom.}, 59(3):353--437, 2001.

\bibitem[HM20]{MiaoHirsch}
S.~Hirsch and P.~Miao.
\newblock A positive mass theorem for manifolds with boundary.
\newblock {\em Pacific J. Math.}, 306(1):185--201, 2020.

\bibitem[HMZ01]{HijaziMontielZhang1}
O.~Hijazi, S.~Montiel, and X.~Zhang.
\newblock Dirac operator on embedded hypersurfaces.
\newblock {\em Math. Res. Lett.}, 8(1-2):195--208, 2001.

\bibitem[HMZ02]{HijaziMontielZhang2}
O.~Hijazi, S.~Montiel, and X.~Zhang.
\newblock Conformal lower bounds for the {D}irac operator of embedded
  hypersurfaces.
\newblock {\em Asian J. Math.}, 6(1):23--36, 2002.

\bibitem[Lic63]{Lichnerowicz2}
A.~Lichnerowicz.
\newblock Spineurs harmoniques.
\newblock {\em C. R. Acad. Sci. Paris}, 257:7--9, 1963.

\bibitem[LL15]{LeeLeFloch}
D.~A. Lee and P.~G. LeFloch.
\newblock The positive mass theorem for manifolds with distributional
  curvature.
\newblock {\em Comm. Math. Phys.}, 339(1):99--120, 2015.

\bibitem[MMT20]{MantoulidisMiaoTam}
C.~Mantoulidis, P.~Miao, and L.-F. Tam.
\newblock Capacity, quasi-local mass, and singular fill-ins.
\newblock {\em J. Reine Angew. Math.}, 768:55--92, 2020.

\bibitem[PT82]{ParkerTaubes}
T.~Parker and C.~H. Taubes.
\newblock On {W}itten's proof of the positive energy theorem.
\newblock {\em Comm. Math. Phys.}, 84(2):223--238, 1982.

\bibitem[Rau08]{Raulot3}
S.~Raulot.
\newblock Rigidity of compact {R}iemannian spin manifolds with boundary.
\newblock {\em Lett. Math. Phys.}, 86(2-3):177--192, 2008.

\bibitem[Rau21]{Raulot12}
S.~Raulot.
\newblock A spinorial proof of the rigidity of the {R}iemannian {S}chwarzschild
  manifold.
\newblock {\em Classical and Quantum Gravity}, 38(8):085015, 2021.

\bibitem[ST02]{ShiTam1}
Y.~Shi and L.-F. Tam.
\newblock Positive mass theorem and the boundary behaviors of compact manifolds
  with nonnegative scalar curvature.
\newblock {\em J. Differential Geom.}, 62(1):79--125, 2002.

\bibitem[SY79a]{SchoenYau2}
R.~M. Schoen and S.-T. Yau.
\newblock Complete manifolds with nonnegative scalar curvature and the positive
  action conjecture in general relativity.
\newblock {\em Proc. Nat. Acad. Sci. U.S.A.}, 76(3):1024--1025, 1979.

\bibitem[SY79b]{SchoenYau1}
R.~M. Schoen and S.-T. Yau.
\newblock On the proof of the positive mass conjecture in general relativity.
\newblock {\em Comm. Math. Phys.}, 65(1):45--76, 1979.

\bibitem[SY92]{SmithYang}
P.~D. Smith and D.~Yang.
\newblock Removing point singularities of {R}iemannian manifolds.
\newblock {\em Trans. Amer. Math. Soc.}, 333(1):203--219, 1992.

\bibitem[Wit81]{Witten1}
E.~Witten.
\newblock A new proof of the positive energy theorem.
\newblock {\em Comm. Math. Phys.}, 80(3):381--402, 1981.

\end{thebibliography}


\end{document}